\documentclass[11 pt]{amsart}
\usepackage{amssymb}
\usepackage{amsmath}
\usepackage{amsfonts}
\usepackage{graphicx}
\usepackage{amsthm}
\usepackage{enumerate}
\usepackage[mathscr]{eucal}
\usepackage{verbatim}

 \theoremstyle{plain}

\theoremstyle{plain}

\theoremstyle{remark}

\def\bbR{{\mathbb {R}}}

\begin{document}

\date{October, 2010}

\title[]
{A uniform Fourier restriction theorem for surfaces in $\bbR  ^d$}

\begin{abstract}
We prove a Fourier restriction result, uniform over a certain collection of reference measures, for 
some indices in the Stein-Tomas range.
\end{abstract}

\author[]
{Daniel M. Oberlin}

\address
{D. M.  Oberlin \\
Department of Mathematics \\ Florida State University \\
Tallahassee, FL 32306}
\email{oberlin@math.fsu.edu}


\maketitle

Let $\sigma$ stand for Lebesgue measure on the unit sphere $S^{d-1}\subset \bbR^d$. The Stein-Tomas Fourier restriction theorem is the estimate 
\begin{equation*}
\|\hat f \|_{L^q (\sigma )}\leq C_p\, \|f\|_p
\end{equation*}
for $f\in L^p (\bbR^d )$ and $1\leq p\leq 2(d+1)/(d+3)$, $q=(d-1)p' /(d+1)$ (where $p$ and $p'$ are conjugate exponents). There is a well-known generalization in which $\sigma$ is replaced by surface area measure on a compact subset with nonvanishing Gaussian curvature of a $(d-1)$-dimensional submanifold 
$\Gamma$ of $\bbR^d$. 
It is also of interest to investigate the possibility of similar Fourier restriction theorems when the Gaussian curvature is allowed to vanish - see \cite{CZ} for a nice discussion of this. If one hopes to get a subset of the Stein-Tomas range of exponents in such a (degenerate)
case, then it follows from results in \cite{IL} that
surface area measure needs to be modified by introducing a weight which vanishes at degenerate
points. When $\Gamma$ is the graph
\begin{equation*}
\{\big(x,\phi (x)\big):x\in\Omega\}
\end{equation*}
of a $C^{(2)}$ function $\phi$ defined on an open subset $\Omega\subset\bbR^{d-1}$, a natural choice of
measure $\lambda$ to replace surface area measure is the so-called affine surface area measure on
$\Gamma$ given by 
\begin{equation*}
d\lambda =|\det H_\phi (x)|^{\frac{1}{d+1}}\,dx .
\end{equation*}
Here is a result for such $\lambda$:

\noindent {\bf Theorem.}  {\sl
Suppose $\Omega$, $\phi$, and $\lambda$ are as above.
Define $\gamma :\Omega \rightarrow \bbR^d$ by $\gamma (x)=\big(x,\phi (x)\big)$. 
There is a positive 
constant $C$ depending only on $d$ and the generic multiplicities of the maps $\text{grad}\, \phi$ and
\begin{equation*}
\Psi :\Omega ^d \rightarrow (\bbR^d )^{d-1},\ \ \Psi (x_1 ,\dots ,x_d )=\big(\gamma (x_1 )-\gamma (x_2 ),
\dots ,\gamma (x_1 )-\gamma (x_d )\big)
\end{equation*}
such that if $q=4(d-1)/(d+1)$, then the Fourier restriction estimate 
\begin{equation}\label{estimate}
\|\hat f \|_{L^{q,\infty}(\lambda )}
\leq C\,
\|f\|_{{4}/{3}}
\end{equation}
holds for $f\in L^{4/3}(\bbR^d )$.}
%

\noindent (We say that a mapping $\Phi$ into $\bbR^n$ has generic multiplicity bounded by $l$ if
$\text{card}\big( \Phi ^{-1}(x)\big)\leq l$ for almost all $x\in\bbR^n$.)

\noindent   {\bf Comments:}

\noindent(a) The papers \cite{O3}, \cite{AS}, and \cite{CKZ} deal with analogues of the Stein-Tomas theorem in three dimensions.
If $d=3$ the conclusion of our result is an $L^{4/3}\rightarrow L^{2,\infty}$
estimate, while the Stein-Tomas theorem gives at the endpoint the stronger $L^{4/3}\rightarrow L^{2}$ estimate.
The range of exponents furnished by interpolating \eqref{estimate}
with the trivial $L^1 \rightarrow L^\infty$ estimate is, for any $d$, a subset of the Stein-Tomas range. 

\noindent (b) Examples in \cite{CZ} show that the multiplicity hypotheses 
are necessary.

\noindent (c) An interesting feature of our result is that the bound it furnishes is uniform 
modulo the multiplicity hypotheses. In particular, B\'ezout's theorem shows that if $n$ is fixed
and if $\phi$ is a polynomial of degree $n$ on $\bbR^{d-1}$, then the constant in the estimate \eqref{estimate} may be chosen independently of $\phi$.

\noindent (d) The (very short) proof below is similar to material in \cite{O2} and \cite{O4}.

\begin{proof}
The constant $C$ may vary from line to line but will depend only on $d$ and the multiplicities mentioned in the hypotheses.  
For the proof we require the estimate
\begin{equation}\label{geometric2}
\int_A \int_A  \chi_E (y_1 -y_2 )\, d\lambda (y_2 )\, d\lambda (y_1 )\leq 
C \, \lambda (A)^{(d-1)/d} \,|E|^{(d-1)/d}
\end{equation}
for Borel sets $A\subset \gamma (\Omega )$, $E\subset \bbR^d$.
Writing $\omega (x)$ for $|\det H_\phi (x)|^{{1}/{d+1}}$, we will deduce \eqref{geometric2} from a
geometric inequality which holds for Borel sets $E\subset \bbR^d$:
\begin{equation}\label{geometric}
\int_{\Omega}\Big[ \int_{ \{ \omega (x_2 )\geq\omega (x_1 )\}}
 \chi_E \big(\gamma (x_1) -\gamma (x_2 )\big)\,\omega (x_2 )\, dx_2 \Big]^d 
\omega (x_1 )\, dx_1
 \leq C\, |E|^{d-1}.
\end{equation}
Estimate \eqref{geometric} is proved in \cite{O1} - see (2) there.
To establish \eqref{geometric2} above, identify a set $A\subset \Omega$ with its image $\gamma (A)$ 
in $\Gamma$ and 
note that the bound 
\begin{multline}\label{ineq0}
\int_A \int_A
\chi_{\{ \omega (x_2 )\geq\omega (x_1 )\}}
\chi_E \big(\gamma (x_1) -\gamma (x_2 )\big)\,\omega (x_2 )\, dx_2 \,
\omega (x_1 )\, dx_1
 \leq   \\
  C\, \lambda (A)^{(d-1)/d} \,|E|^{(d-1)/d}
\end{multline}
follows from \eqref{geometric}. Replacing $E$ by $-E$ and interchanging $x_1$ and $x_2$ gives 
\begin{multline*}
\int_A \int_A
\chi_{\{ \omega (x_2 )\leq\omega (x_1 )\}}
\chi_E \big(\gamma (x_1) -\gamma (x_2 )\big)\,\omega (x_2 )\, dx_2 \,
\omega (x_1 )\, dx_1
 \leq   \\
  C\, \lambda (A)^{(d-1)/d} \,|E|^{(d-1)/d}.
\end{multline*}
With \eqref{ineq0} this implies \eqref{geometric2}. For a measure $\mu$ on $\bbR^d$, we will write $\widetilde{\mu}$ for the measure defined by $\widetilde{\mu}(E)=\mu (-E)$ and will then interpret \eqref{geometric2} as the 
Lorentz space estimate 
\begin{equation}\label{ineq1}
\|(\chi_A d\lambda )\ast \widetilde{(\chi_A d\lambda )}\|_{d,\infty}\leq C\,\lambda (A)^{(d-1)/d} .
\end{equation}
Interpolating \eqref{ineq1} with 
\begin{equation*}
\|(\chi_A d\lambda )\ast \widetilde{(\chi_A d\lambda )}\|_{1}= \lambda (A)^2
\end{equation*}
yields 
\begin{equation*}
\|\, |\widehat{\chi_A d\lambda}|^2 \|_2 =
\|(\chi_A d\lambda )\ast \widetilde{(\chi_A d\lambda )}\|_{2}\leq C\,\lambda (A)^{(3d-5)/(2d-2)}
\end{equation*}
and so 
\begin{equation*}
\|\widehat{\chi_A d\lambda}\|_4 \leq C\,\lambda (A)^{(3d-5)/(4d-4)}.
\end{equation*}
This is the dual of \eqref{estimate}.

\end{proof}

\end{document}